\documentclass[smallextended]{svjour3}
\usepackage{amsmath}
\usepackage{graphicx}
\usepackage{array,a4,amssymb}
\newtheorem{aaaaa}{DO NOT USE!}
\newtheorem{Corollary}[aaaaa]{Corollary}

\newtheorem{Definition}[aaaaa]{Definition}
\newtheorem{Example}[aaaaa]{Example}

\newtheorem{Proposition}[aaaaa]{Proposition}
\newtheorem{Remark}[aaaaa]{Remark}
\newtheorem{Theorem}[aaaaa]{Theorem}
\usepackage{graphicx,booktabs,url,color}
\usepackage{multirow}
\newcommand{\startprf}{\noindent\emph{Proof:} }
\providecommand{\qed}{\hfill\rule{1.5ex}{1.5ex}}
\newcommand{\startprfoftheorem}{\noindent\emph{Proof of Theorem~\ref{prosmall}}: }

\begin{document}
\begin{center}
{\Large{More results on the number of zeros of multiplicity at least $r$}}\\
\ \\
Olav Geil and Casper Thomsen\\
Department of Mathematical Sciences\\
Aalborg University\\
Fr.\ Bajersvej 7G\\
9220 Aalborg {\O}\\
Denmark\\
\ \\

Email: olav@math.aau.dk and ct@spag.dk
\end{center}
\ \\
\begin{center}
\begin{minipage}{12cm}
{\small{{\textbf{Abstract:}}
    \     We consider multivariate polynomials and investigate how many zeros of multiplicity at least $r$ they can have over a Cartesian product of finite subsets of a field. Here $r$ is any prescribed positive integer and the 
    definition of
    multiplicity that we use is the one related to Hasse derivatives. As a generalization of material in~\cite{augot,dvir} a general version of the Schwartz-Zippel was presented in~\cite{weighted} which from the leading monomial -- with respect to a lexicographic ordering -- estimates the sum of zeros when counted with multiplicity. The corresponding corollary on the number of zeros of multiplicity at least $r$ is in general not sharp and therefore in~\cite{weighted} a recursively defined function $D$ was introduced using which one can derive improved information. The recursive function being rather complicated, the only known closed formula consequences of it are for the case of two variables~\cite{weighted}. In the present paper we derive closed formula consequences for arbitrary many variables, but for the powers in the leading monomial being not too large.  Our bound
    can be viewed as a generalization of the footprint bound~\cite{onorin,gh} -- the classical footprint bound taking not multiplicity into account.\\
    
\ \\

\noindent {\textbf{Keywords:}} Footprint bound, multiplicity, multivariate polynomial,
Schwartz-Zippel bound, zeros of polynomial \\

\noindent {\textbf{MSC classifications:}} Primary: 12Y05. Secondary:
11T06, 12E05, 13P05, 26C99 }}

\end{minipage}
\end{center}
\ \\

\section{Introduction}\label{sec1} 
Given a univariate polynomial over an arbitrary field it is an easy
task to estimate the number of zeros of multiplicity at least $r$, for any fixed positive integer $r$. As is well-known the number of such zeros is less than or equal to the
degree of the polynomial divided by $r$. For multivariate
polynomials the situation is much more complicated as these polynomials
on the one hand typically have an infinite number of zeros when the
field is infinite and on the other hand have only a finite number of
zeros when not. A meaningful reformulation of the problem which works
independently of the field -- and which will be taken in the present
paper -- is to restrict to point sets that are Cartesian products of
finite sets. This of course includes the important case where the
point set is ${\mathbb{F}}_q \times \cdots \times {\mathbb{F}}_q$,
${\mathbb{F}}_q$ being the finite field with $q$ elements. Another
concern is which definition of multiplicity to use as for multivariate
polynomials there are more competing definitions. In the present paper
we use the one related to Hasse derivatives (see
Definition~\ref{defmult} below).\\

 The interest in studying the outlined
problem originally came from applications to Guruswami-Sudan style \cite{GS} list decoding algorithms for
$q$-ary Reed-Muller codes, weighted Reed-Muller codes and their likes
\cite{pw_ieee,pw_expanded,manyauthors,augot,weighted}. The first bound on the number of zeros of prescribed
multiplicity was developed by Pellikaan and Wu in~\cite{pw_ieee,pw_expanded}. Later Augot
and Stepanov improved upon Pellikaan and Wu's bound (see~\cite[Prop.\ 13]{weighted}) by generalizing the Schwartz-Zippel bound to also deal with
multiplicity \cite{augot}. The proof of this bound was later given by Dvir et
al.\ in~\cite{dvir} where it was used to estimate the size of 
Kakeya sets over finite fields. The mentioned Schwartz-Zippel bound estimates the sum of zeros when counted with multiplicity. From this, one obtains an easy corollary on the number of zeros of multiplicity $r$ or more. All of the above mentioned bounds  are  
stated in terms of the total degree of the involved polynomials
and the point set under consideration is always ${\mathbb{F}}_q \times
\cdots \times 
{\mathbb{F}}_q$. In~\cite[Th.\ 5]{weighted} the generalization of the Schwartz-Zippel
bound was taken a step further to now work for arbitrary finite point
sets $S_1 \times \cdots \times S_m$, $S_i \subseteq {\mathbb{F}}$, $i=1, \ldots , m$
(where ${\mathbb{F}}$ is any field) and to take into account the
leading monomial with respect to a lexicographic ordering. Again one
obtains an easy corollary on the number of zeros of multiplicity at
least $r$ \cite[Cor.\ 3]{weighted}. Whereas the generalized
Schwartz-Zippel bound~\cite[Th.\ 5]{weighted} is tight in the sense
that we can always find polynomials attaining it (see
Proposition~\ref{prosharp} below) a similar result does not hold for
its corollary~\cite[Cor.\ 3]{weighted}. To address this problem we
introduced in~\cite{weighted} a recursively defined function $D$ to
estimate the number of zeros of multiplicity at least
$r$. Unfortunately, the function $D$ is quite complicated and only for
the case of two variables some simple closed formula upper bounds were derived~\cite[Prop.\ 16]{weighted}.

The purpose of the present paper is to establish for the general case of arbitrarily many variables a class of cases in which from $D$ we can derive a simple closed formula expression which is still an improvement to the Schwartz-Zippel bound for zeros of multiplicity at least $r$ (\cite[Cor.\ 3]{weighted}). The bound that we derive turns out to be a natural generalization of the
footprint bound \cite{onorin,gh} which estimates the number of zeros without
taking multiplicity into consideration.\\

The paper is organized as follows. In Section~\ref{sec2} we start by defining multiplicity and by recalling the general Schwartz-Zippel bound and as a corollary the Schwartz-Zippel bound for zeros of multiplicity at least $r$. The rest of Section~\ref{sec2} is devoted to a discussion of the method from~\cite{weighted}. In Section~\ref{sec3} we give the new results regarding a simple closed formula upper bound for the case of the coefficients in the leading monomial being small. The concept of being small in general is rather involved and we therefore establish simple sufficient conditions for this to happen.

\section{Background}\label{sec2}

We first recall the concept of Hasse derivatives.

\begin{Definition}
Given $F(X_1, \ldots  X_m)\in {\mathbb{F}}[X_1, \ldots , X_m]$ and
$\vec{k}=(k_1, \ldots , k_m) \in {\mathbb{N}}_0^m$ the $\vec{k}$'th
Hasse derivative of $F$, denoted by $F^{(\vec{k})}(X_1, \ldots , X_m)$ is the
coefficient of $Z_1^{k_1} \cdots Z_m^{k_m}$ in $F(X_1+Z_1, \ldots , X_m+Z_m)\in {\mathbb{F}}(X_1, \ldots , X_m)[Z_1, \ldots , Z_m]$. In other words 
$$F(X_1+Z_1, \ldots , X_m+Z_m)=\sum_{\vec{k}} F^{(\vec{k})}(X_1, \ldots , X_m)Z_1^{k_1}
\cdots Z_m^{k_m}.$$
\end{Definition}

Observe that the next definition includes the usual concept of multiplicity for univariate polynomials as a special case.

\begin{Definition}\label{defmult}
For $F(X_1, \ldots , X_m) \in {\mathbb{F}}[X_1, \ldots , X_m]\backslash \{ {0} \}$ and
$\vec{a}=(a_1, \ldots , a_m) \in {\mathbb{F}}^m$ we define the multiplicity of $F$ at $\vec{a}$
denoted by ${\mbox{mult}}(F,\vec{a})$ as follows. Let $r$ be an
integer such that for every $\vec{k}=(k_1, \ldots ,
k_m) \in {\mathbb{N}}_0^m$ with $k_1+\cdots +k_m < r$, $F^{(\vec{k})}(a_1, \ldots , a_m)=0$
  holds, but for some  $\vec{k}=(k_1, \ldots ,
k_m) \in {\mathbb{N}}_0^m$ with $k_1+\cdots +k_m = r$,
$F^{(\vec{k})}(a_1, \ldots , a_m)\neq 0$ holds, then 
${\mbox{mult}}(F,\vec{a})=r$. If $F=0$ then we define ${\mbox{mult}}(F,\vec{a})=\infty$.
\end{Definition}

The above definition is the one that is usually given in the literature. For our purpose the below equivalent description shall also  prove useful.

\begin{Definition}\label{defeq}
Let $F(X_1, \ldots , X_m) \in {\mathbb{F}}[X_1, \ldots , X_m]\backslash \{0\}$ and
$\vec{a}=(a_1, \ldots , a_m) \in {\mathbb{F}}^m$. 
Consider the ideal 
\begin{eqnarray}
J_t=\langle (X_1-a_1)^{p_1}\cdots(X_m-a_m)^{p_m} \mid p_1+\cdots
+p_m=t \rangle \subseteq {\mathbb{F}}[X_1, \ldots , X_m]. \nonumber
\end{eqnarray}
We have
${\mbox{mult}}(F,\vec{a})=r$ if $F\in J_r \backslash
J_{r+1}$. If $F=0$ we have ${\mbox{mult}}(F,\vec{a})=\infty$.
\end{Definition}

We next state the most general form of the Schwartz-Zippel bound for fields
\cite[Th.\ 5]{weighted}. Here, and in the rest of the paper $S_1,
\ldots , S_m \subset {\mathbb{F}}$ are finite subsets of the field
${\mathbb{F}}$ and we write $s_1 =|S_1|, \ldots , s_m=|S_m|$. We note
that the below theorem was generalized to arbitrary commutative rings
in~\cite[Th.\ 7.10]{bishnoi2015zeros} where it was called the
generalized Schwartz Theorem.
\begin{Theorem}\label{prop-sz-gen}
Let $F(X_1, \ldots , X_m) \in {\mathbb{F}}[X_1, \ldots , X_m]$ be a non-zero polynomial and
let $X_1^{i_1} \cdots X_m^{i_m}$ be its leading
monomial with respect to a lexicographic ordering $\prec_{lex}$. Then for
any finite sets $S_1, \ldots ,S_m \subseteq {\mathbb{F}}$
\begin{eqnarray}
\sum_{\vec{a}\in S_1 \times \cdots \times S_m}{\mbox{mult}}(F,\vec{a})
\leq i_1s_2\cdots s_m+s_1i_2s_3 \cdots s_m+\cdots +s_1\cdots s_{m-1}i_m.\nonumber
\end{eqnarray}
\end{Theorem}

Turning to the problem of estimating the number of zeros of multiplicity at least $r$ -- which is the topic of the present paper -- we have the following corollary corresponding to~\cite[Cor.\ 3]{weighted}. We may think of it as the Schwartz-Zippel bound for zeros of multiplicity at least $r$. 
\begin{Corollary}\label{cor-sz-gen}
Let $F(X_1, \ldots , X_m) \in {\mathbb{F}}[X_1, \ldots , X_m]$ be a non-zero polynomial and
let $X_1^{i_1} \cdots X_m^{i_m}$ be its leading
monomial with respect to the lexicographic ordering. Assume $S_1, \ldots
,S_m \subseteq {\mathbb{F}}$ are finite sets. 
Then over
$S_1 \times \cdots \times S_m $ the number
of zeros of multiplicity at least $r$ is less than or equal to the
minimum of 
\begin{eqnarray}
\big( i_1s_2\cdots s_m+s_1i_2s_3\cdots s_m+\cdots +s_1\cdots
s_{m-1}i_m \big) /r\nonumber
\end{eqnarray}
and $s_1 \cdots s_m$.
\end{Corollary}

As mentioned in the introduction one 
obtains better estimates than Corollary~\ref{cor-sz-gen}
by using the recursively defined function $D$. In particular
Corollary~\ref{cor-sz-gen} is not tight. Before giving the details we
pause for a moment to show that on the other hand
Theorem~\ref{prop-sz-gen} is tight (a fact that has not been reported
before).  For this purpose we shall need the notation 
$$S_j=\{\alpha_1^{(j)}, \ldots , \alpha_{s_j}^{(j)} \}$$
for $j=1, \ldots , m$, and the below proposition:
\begin{Proposition}\label{pronuogsaa}
Consider 
\begin{equation}
  F(X_1, \ldots , X_m)=\prod_{u=1}^{m} \prod_{v=1}^{s_u}(X_u-\alpha_v^{(u)})^{r_{v}^{(u)}}.\label{eqlinfac}
\end{equation}
The multiplicity of $(\alpha_{j_1}^{(1)}, \ldots ,
\alpha_{j_m}^{(m)})$ in $F(X_1, \ldots , X_m)$ equals
\begin{equation}
r_{j_1}^{(1)}+\cdots +r_{j_m}^{(m)}. \label{eqHjubi}
\end{equation}
\end{Proposition}
\startprf
Clearly, the multiplicity is greater than or equal to
$r=r_{j_1}^{(1)}+\cdots +r_{j_m}^{(m)}$. Using Gr\"{o}bner basis theory we now
show that it is not larger. We substitute
${\mathcal{X}}_i=X_i-\alpha_{j_i}^{(i)}$ for $i=1, \ldots ,m$ and observe
that by Buchberger's S-pair criteria 
$${\mathcal{B}}=\{{\mathcal{X}}_1^{r_1} \cdots {\mathcal{X}}_m^{r_m}
\mid r_1+\cdots +r_m=r+1\}$$
is a Gr\"{o}bner basis (with respect to any fixed monomial
ordering). 
The support of $F({\mathcal{X}}_1, \ldots , {\mathcal{X}}_m)$ contains
a monomial of the form ${\mathcal{X}}_1^{i_1} \cdots
{\mathcal{X}}_m^{i_m}$ with $i_1+\cdots +i_m=r$. 
 Therefore the remainder of $F({\mathcal{X}}_1, \ldots ,
{\mathcal{X}}_m)$ modulo ${\mathcal{B}}$ is non-zero. 
It is well known that if a
polynomial is reduced modulo a Gr\"{o}bner basis then the remainder is
zero if and only if it belongs to the ideal generated by the elements
in the basis. \qed
\ \\

We are now ready to show that Theorem~\ref{prop-sz-gen} is tight. 
\begin{Proposition}\label{prosharp}
Let $S_1, \ldots , S_m \subseteq {\mathbb{F}}$ be finite sets. If
$F(X_1, \ldots , X_m)\in {\mathbb{F}}[X_1, \ldots , X_m]$ is a product
of univariate linear factors -- meaning that it is of the form (\ref{eqlinfac}) -- then the number of
zeros of $F$ counted with multiplicity reaches the generalized Schwartz-Zippel bound
(Theorem~\ref{prop-sz-gen}).
\end{Proposition}
\startprf
Consider the polynomial
$$F(X_1, \ldots , X_m)=\prod_{u=1}^m \prod_{v=1}^{s_u}\big(X_u-\alpha_v^{(u)}\big)^{r_v^{(u)}}.$$
Write $i_u=\sum_{v=1}^{s_u}r_v^{(u)}$, $u=1, \ldots , m$. Applying carefully Proposition~\ref{pronuogsaa} we obtain
\begin{eqnarray}
\sum_{\vec{a} \in S_1 \times \cdots \times S_m}
{\mbox{mult}}(F,\vec{a})
&=\sum_{t=1}^{s_1}(s_2 \cdots s_m  )r_t^{(1)}+\cdots +
\sum_{t=1}^{s_m}(s_1 \cdots s_{m-1}  )r_t^{(m)} \nonumber
\\
&=i_1s_2  \cdots s_{m}  +\cdots +s_1 
\cdots s_{m-1} i_m \nonumber
\end{eqnarray}
and we are through. \qed

\ \\
We next return to the problem of improving 
Corollary~\ref{cor-sz-gen} for which we introduced
in~\cite[Def.\ 5]{weighted}
the function $D$.

\begin{Definition}\label{defD}
Let $r \in {\mathbb{N}}, i_1, \ldots , i_m \in {\mathbb{N}}_0$. Define 
$$D(i_1,r,s_1)=\min \big\{\big\lfloor \frac{i_1}{r} \big\rfloor,s_1\big\}$$
and for $m \geq 2$
\begin{multline*}
D(i_1, \ldots , i_m,r,s_1, \ldots ,s_m)=
\\
\begin{split}
\max_{(u_1, \ldots  ,u_r)\in A(i_m,r,s_m) }&\bigg\{ (s_m-u_1-\cdots -u_r)D(i_1,\ldots ,i_{m-1},r,s_1,
\ldots ,s_{m-1})\\
&\quad+u_1D(i_1, \ldots , i_{m-1},r-1,s_1, \ldots ,s_{m-1})+\cdots
\\
&\quad +u_{r-1}D(i_1, \ldots ,i_{m-1},1,s_1, \ldots , s_{m-1})+u_rs_1\cdots
s_{m-1} \bigg\}
\end{split}
\end{multline*}
where 
\begin{multline}
A(i_m,r,s_m)= \\
\{ (u_1, \ldots , u_r) \in {\mathbb{N}}_0^r \mid u_1+ \cdots
+u_r \leq s_m {\mbox{ \ and \ }} u_1+2u_2+\cdots +ru_r \leq i_m\}.\label{eqA}
\end{multline}
\end{Definition}
Throughout the rest of the paper we shall always assume that $r \in
{\mathbb{N}}$ and that $i_1, \ldots , i_m \in {\mathbb{N}}_0$. 
The improvement of Corollary~\ref{cor-sz-gen} was given in~\cite[Th.\ 6]{weighted}
as follows:

\begin{Theorem}\label{prorec}
For a polynomial $F(X_1, \ldots , X_m)\in {\mathbb{F}}[X_1, \ldots , X_m]$ let $X_1^{i_1}\cdots X_m^{i_m}$ be its leading monomial with
respect to the lexicographic ordering $\prec_{lex}$ with $X_m  \prec_{lex} \cdots \prec_{lex} X_1$. Then $F$ has at most $D(i_1, \ldots , i_m,r,s_1,
\ldots ,s_m)$ zeros of multiplicity at least $r$ in $S_1\times \cdots
\times S_m$. The corresponding recursive algorithm produces a number
that is at most equal to the number found in
Corollary~\ref{cor-sz-gen} and is at most equal to $s_1 \cdots s_m$.
\end{Theorem}

When
$\lfloor i_1/s_1 \rfloor + \cdots + \lfloor i_m/s_m \rfloor \geq r$
Proposition~\ref{pronuogsaa} guarantees the existence of polynomials $F(X_1, \ldots , X_m)$ with leading monomial $X_1^{i_1} \cdots X_m^{i_m}$ having all elements of $S_1 \times \cdots \times S_m$ as zeros of multiplicity at least $r$. Hence, we only need to apply Theorem~\ref{prorec} to the case
$
\lfloor i_1/s_1 \rfloor + \cdots + \lfloor i_m/s_m \rfloor < r, $
and in particular we can assume $i_t < r s_t$.

\begin{Example}\label{exny1}
In this example we estimate the number of zeros of multiplicity 
$3$ or more for polynomials in two
variables. Both $S_1$ and $S_2$ are assumed to be of size $5$. 
From the above discussion, for 
\begin{eqnarray}
(i_1,i_2)&\in&\{(\alpha,\beta) \mid \alpha \geq 15\} \cup \{ (\alpha,\beta) \mid \alpha \geq 10
  {\mbox{ and }} \beta \geq 5\}\nonumber \\
&&\cup \{(\alpha,\beta) \mid \alpha \geq 5 {\mbox{ and }} \beta \geq 10\} \cup \{(\alpha,\beta)
  \mid \beta \geq 15\}\nonumber 
\end{eqnarray}
we have $D(i_1,i_2,3,5,5)=25$. Table~\ref{tabny1} shows 
information obtained from our algorithm for the remaining possible
choices of exponents $(i_1,i_2)$. Observe, that the table is not symmetric
meaning that $D(i_1,i_2,3,5,5)$ does not always equal
$D(i_2,i_1,3,5,5)$. The corresponding values of the
Schwartz-Zippel bound 
(Corollary~\ref{cor-sz-gen}) is displayed in  Table~\ref{tabny2}, from
which it is clear that indeed the function $D$ can sometimes give a
dramatic improvement. For instance $D(3,11,3,5,5)$ equals $19$, but the
Schwartz-Zippel bound only gives the estimate $23$. Similarly,
$D(2,12)$ equals $20$ and the Schwartz-Zippel bound gives $23$.
\begin{table}
\centering
\caption{$D(i_1,i_2,3,5,5)$}
\newcommand{\SP}{~~}
\begin{tabular}{@{}c@{}cc@{\SP}c@{\SP}c@{\SP}c@{\SP}c@{\SP}c@{\SP}c@{\SP}c@{\SP}c@{\SP}c@{\SP}c@{\SP}c@{\SP}c@{\SP}c@{\SP}c@{}}
    \toprule
    &&\multicolumn{15}{c}{$i_1$}\\
    &&0   &1   &2   &3   &4   &5   &6   &7   &8   &9   &10  &11  &12  &13  &14\\
    \addlinespace
    \multirow{15}{*}{$i_2$}
    &\multicolumn{1}{c}{0 }&0 &0 &0 &5 &5 &5 &10&10&10&15&15&15&20&20&20\\
    &\multicolumn{1}{c}{1 }&0 &0 &1 &5 &6 &6 &11&11&12&{{16}}&17&17&{{21}}&21&21\\
    &\multicolumn{1}{c}{2 }&0 &1 &2 &7 &8 &9 &{{13}}&13&14&{{17}}&{{19}}&19&{{22}}&22&22\\
    &\multicolumn{1}{c}{3 }&5 &5 &5 &9 &9 &10&{{14}}&14&{{16}}&{{18}}&{21}&{21}&{23}&{23}&23\\
    &\multicolumn{1}{c}{4 }&5 &5 &6 &9 &11&{13}&{16}&{16}&{18}&{19}&{23}&{23}&{24}&{24}&{24}\\
    &\multicolumn{1}{c}{5 }&5 &6 &7 &11&12&{14}&{17}&{17}&{20}&{20}\\
    &\multicolumn{1}{c}{6 }&10&10&10&13&14&{17}&{19}&{19}&{21}&{21}\\
    &\multicolumn{1}{c}{7 }&10&10&11&13&15&{18}&{20}&{20}&{22}&{22}\\
    &\multicolumn{1}{c}{8 }&10&11&12&15&{17}&{21}&{22}&{22}&{23}&{23}\\
    &\multicolumn{1}{c}{9 }&15&15&15&17&{18}&{22}&{23}&{23}&{24}&{24}\\
    &\multicolumn{1}{c}{10}&15&15&16&17&{20}\\
    &\multicolumn{1}{c}{11}&15&16&17&19&{21}\\
    &\multicolumn{1}{c}{12}&20&20&20&21&{22}\\
    &\multicolumn{1}{c}{13}&20&20&21&21&{23}\\
    &\multicolumn{1}{c}{14}&20&21&22&23&{24}\\
    \bottomrule
\end{tabular}
\label{tabny1}
\end{table}

\begin{table}
\centering
\caption{The Schwartz-Zippel bound (sz) for zeros of multiplicity at least $3$}
\begin{tabular}{c|rrrrrrrrrrrr}
$i_1+i_2$&0&1&2&3&4&5&6&7&8&9&10&11\\
sz&0   &1   &3   &5   &6   &8   &10   &11   &13   &15   &16  &18\\
\ \\
$i_1+i_2$&12&13&14&15&16&17&18\\
 sz &20
  &21  &23&25&25&25&25
\end{tabular}
\label{tabny2}
\end{table}

It is easy to  establish a lower bound on the maximal number of
possible zeros of multiplicity at least $r=3$ for polynomials with any leading monomial
$X_1^{i_1}X_2^{i_2}$. This is done by inspecting polynomials of the
form~(\ref{eqlinfac}). As an example 
$\prod_{u=1}^{4}(X_1-\alpha_{u}^{(1)})^2\prod_{v=1}^5(X_2-\alpha_v^{(2)})$
has $20$ zeros of multiplicity (at least) $3$. But $D(8,5,3,5,5)=20$ and
therefore the true value of the maximal number of zeros of
multiplicity at least $3$ is $20$ in this case. In Table~\ref{tabnyyy3} we
list the difference between $D(i_1,i_2,3,5,5)$ and the lower bound
found by using the above method. The large amount of zero's in the table proves that
$D(i_1,i_2,3,5,5)$ often equals the true maximal number of zeros of
multiplicity at least $3$.
\begin{table}
\centering
\caption{Difference between upper and lower bound in Example~\ref{exny1}}\label{tabny6}
\newcommand{\SP}{~\hspace*{0.851em}}
\begin{tabular}{@{}c@{}cc@{\SP}c@{\SP}c@{\SP}c@{\SP}c@{\SP}c@{\SP}c@{\SP}c@{\SP}c@{\SP}c@{~~~}c@{~~}c@{~~}c@{~~}c@{~~}c@{}}
    \toprule
    &&\multicolumn{15}{c}{$i_1$}\\
    &&0&1&2&3&4&5&6&7&8&9&10&11&12&13&14\\
    \addlinespace
    \multirow{15}{*}{$i_2$}
    &\multicolumn{1}{c}{0} &0&0&0&0&0&0&0&0&0&0&0&0&0&0&0\\
    &\multicolumn{1}{c}{1} &0&0&0&0&1&0&1&1&1&1&2&1&1&1&0\\
    &\multicolumn{1}{c}{2} &0&0&0&2&2&2&3&2&2&2&3&2&2&1&0\\
    &\multicolumn{1}{c}{3} &0&0&0&0&0&1&1&1&3&1&4&3&2&2&0\\
    &\multicolumn{1}{c}{4} &0&0&0&0&2&3&3&3&2&2&3&2&2&1&0\\
    &\multicolumn{1}{c}{5} &0&0&0&2&2&3&2&2&0&0\\
    &\multicolumn{1}{c}{6} &0&0&0&0&1&2&3&2&1&0\\
    &\multicolumn{1}{c}{7} &0&0&0&0&2&3&3&3&1&0\\
    &\multicolumn{1}{c}{8} &0&0&0&2&1&1&2&1&2&0\\
    &\multicolumn{1}{c}{9} &0&0&0&0&1&2&2&1&1&0\\
    &\multicolumn{1}{c}{10}&0&0&0&0&0\\
    &\multicolumn{1}{c}{11}&0&0&0&1&0\\
    &\multicolumn{1}{c}{12}&0&0&0&0&0\\
    &\multicolumn{1}{c}{13}&0&0&0&0&0\\
    &\multicolumn{1}{c}{14}&0&0&0&0&0\\
    \bottomrule
\end{tabular}
\label{tabnyyy3}
\end{table}

\end{Example}
In~\cite[Pro.\ 16]{weighted} we derived the following
closed formula expression upper bounds for
the case of two variables. 

\begin{Proposition}\label{protwovar}
For $k=1, \ldots , r-1$,  $D(i_1,i_2,r,s_1,s_2)$ is upper bounded by\\
$\begin{array}{cl}
{\mbox{(C.1)}}&  {\displaystyle{s_2\frac{i_1}{r}+\frac{i_2}{r}\frac{i_1}{r-k}}}\\
&{\mbox{if \  }}(r-k)\frac{r}{r+1}s_1 \leq i_1 < (r-k)s_1
{\mbox{ \ and \ }} 0\leq i_2 <ks_2\\
{\mbox{(C.2)}}&
  {\displaystyle{s_2\frac{i_1}{r}+((k+1)s_2-i_2)(\frac{i_1}{r-k}-\frac{i_1}{r})+(i_2-ks_2)(s_1-\frac{i_1}{r})}}\\
& {\mbox{if \ }}(r-k)\frac{r}{r+1}s_1 \leq i_1 < (r-k)s_1 {\mbox{ \
    and \ }} ks_2\leq i_2 <(k+1)s_2\\
{\mbox{(C.3)}}&
{\displaystyle{s_2\frac{i_1}{r}+\frac{i_2}{k+1}(s_1-\frac{i_1}{r})}}\\
&{\mbox{if \ }} (r-k-1)s_1 \leq i_1 < (r-k)\frac{r}{r+1}s_1 {\mbox{ \
    and \ }} 0 \leq i_2 < (k+1)s_2.
\end{array}
$\\
Finally,\\
$\begin{array}{cl}
{\mbox{(C.4)}}& {\displaystyle{D(i_1,i_2,r,s_1,s_2)=s_2\lfloor \frac{i_1}{r} \rfloor
  +i_2(s_1-\lfloor \frac{i_1}{r} \rfloor )}}\\
& {\mbox{if \ }} s_1(r-1) \leq i_1 < s_1r {\mbox{ \ and \ }} 0 \leq i_2 < s_2.
\end{array}
$\\
The above numbers are at most equal to $\min\{(i_1s_2+s_1i_2)/r, s_1s_2 \}$.
\end{Proposition}
If in (C.3) of the above proposition we substitute $k=r-1$ then we derive
\begin{equation}
  D(i_1,i_2,r,s_1,s_2) \leq s_1s_2-(s_1-\frac{i_1}{r})(s_2-\frac{i_2}{r}) \label{eqmer2}
\end{equation}
for $0\leq i_1<  \frac{r}{r+1}s_1$ and $0 \leq i_2 < r s_2$. Actually, (\ref{eqmer2}) holds under the weaker assumption
\begin{equation}
  0\leq i_1 \leq  \frac{r}{r+1}s_1, 0 \leq i_2 < rs_2 \label{eqbett}
\end{equation}
which is seen by plugging in the values $k=r-1$ and
$i_1=\frac{r}{r+1}s_1$ into the expressions in (C.1), (C.2) and
(\ref{eqmer2}). This is the result that we will generalize to more
variables in the next section.

\begin{Example}\label{extwo}
This is a continuation of Example~\ref{exny1} where we investigated  
$D(i_1,i_2,3,5,5)$. Although
condition~(\ref{eqbett}) reads $i_1 \leq 3$ and $i_2 \leq 14$ we print in
Table~\ref{tabny3} the value of~(\ref{eqmer2}) for all possible
$(i_1,i_2)$. The single, as well as double, underlined numbers correspond to entries where the
number is strictly smaller than 
$D(i_1,i_2,3,5,5)$. For such entries (\ref{eqmer2}) certainly doesn't hold
true. By inspection, condition~(\ref{eqbett}) seems rather sharp. The
double underlined numbers correspond to cases where even,  the value is 
 smaller than the lower bounds on the maximal number of zeros, that
 we established at the end of Example~\ref{exny1}. Hence, not only
 cannot~(\ref{eqmer2}) serve as a general upper bound on $D$, but neither
 can it serve as a general upper bound on the maximal number of zeros
 of multiplicity at least $r$.
\begin{table}
\centering
\caption{$\lfloor 25-(5-i_1/3)(5-i_2/3)\rfloor$}
\newcommand{\SP}{~~}
\begin{tabular}{@{}c@{}cc@{\SP}c@{\SP}c@{\SP}c@{\SP}c@{\SP}c@{\SP}c@{\SP}c@{\SP}c@{\SP}c@{\SP}c@{\SP}c@{\SP}c@{\SP}c@{\SP}c@{}}
    \toprule
    &&\multicolumn{15}{c}{$i_1$}\\
    &&0   &1   &2   &3   &4   &5   &6   &7   &8   &9   &10  &11  &12  &13  &14\\
    \addlinespace
    \multirow{15}{*}{$i_2$}
    &\multicolumn{1}{c}{0 }&0 &1 &3 &5 &\it{6} &\it{8}&\it{10}&\it{11}&\it{13}&\it{15}&\it{16}&\it{18}&\it{20}&\it{21}&\it{23}\\
    &\multicolumn{1}{c}{1 }&1 &3 &4 &6 &\it{7} &\it{9} &\it{11}&\it{12}&{\it{14}}&{\bf{\underline{\it{15}}}}&\it{17}&\it{18}&{\underline{\it{20}}}&\it{21}&\it{23}\\
    &\multicolumn{1}{c}{2 }&3 &4 &6 &7 &\it{9} &\it{10} &{\underline{\it{12}}}&\it{13}&\it{14}&{\underline{\it{16}}}&{\underline{\it{17}}}&\it{19}&{\underline{\it{20}}}&\it{22}&\it{23}\\
    &\multicolumn{1}{c}{3 }&5 &6 &7 &9 &\it{10} &\it{11}&{\underline{\it{13}}}&\it{14}&{\underline{\it{15}}}&{\underline{\it{17}}}&\underline{\it{18}}&\underline{\it{19}}&\underline{\it{21}}&\underline{\it{22}}&\it{23}\\
    &\multicolumn{1}{c}{4 }&6 &7 &9 &10 &\it{11}&\underline{\it{12}}&\underline{\it{14}}&\underline{\it{15}}&\underline{\it{16}}&\underline{\it{17}}&\underline{\underline{\it{18}}}&\underline{\underline{\it{20}}}&\underline{\underline{\it{21}}}&\underline{\underline{\it{22}}}&\underline{\underline{\it{23}}}\\
    &\multicolumn{1}{c}{5 }&8 &9 &10 &11&\it{12}&\underline{\it{13}}&\underline{\it{15}}&\underline{\it{16}}&\underline{\underline{\it{17}}}&\underline{\underline{\it{18}}}\\
    &\multicolumn{1}{c}{6 }&10&11&12&13&\it{14}&\underline{\it{15}}&\underline{\it{16}}&\underline{\it{17}}&\underline{\underline{\it{18}}}&\underline{\underline{\it{19}}}\\
    &\multicolumn{1}{c}{7 }&11&12&13&14&\it{15}&\underline{\it{16}}&\underline{\it{17}}&\underline{\it{17}}&\underline{\underline{\it{18}}}&\underline{\underline{\it{19}}}\\
    &\multicolumn{1}{c}{8 }&13&14&14&15&\underline{\it{16}}&
\underline{\underline{\it{17}}}&\underline{\underline{\it{18}}}&\underline{\underline{\it{18}}}&\underline{\underline{\it{19}}}&\underline{\underline{\it{20}}}\\
    &\multicolumn{1}{c}{9 }&15&15&16&17&\underline{\it{17}}&\underline{\underline{\it{18}}}&\underline{\underline{\it{19}}}&\underline{\underline{\it{19}}}&\underline{\underline{\it{20}}}&\underline{\underline{\it{21}}}\\
    &\multicolumn{1}{c}{10}&16&17&17&18&\underline{\underline{\it{18}}}\\
    &\multicolumn{1}{c}{11}&18&18&19&19&\underline{\underline{\it{20}}}\\
    &\multicolumn{1}{c}{12}&20&20&20&21&\underline{\underline{\it{21}}}\\
    &\multicolumn{1}{c}{13}&21&21&22&22&\underline{\underline{\it{22}}}\\
    &\multicolumn{1}{c}{14}&23&23&23&23&\underline{\underline{\it{23}}}\\
    \bottomrule
\end{tabular}
\label{tabny3}
\end{table}

\end{Example}

{\section{A closed formula expression when  $(i_1, \ldots ,i_m)$ is small}\label{secismall}\label{sec3}}
\noindent

Having already four different cases of closed formula expressions when
$m=2$ (Proposition~\ref{protwovar}), 
the situation gets
very complicated for more variables. Assuming, however, that
the exponent $(i_1, \ldots , i_m)$ in the leading monomial is ``small''
 -- a concept that will be formally defined in Definition~\ref{defconda} below -- we
can give a simple formula which is a generalization of (\ref{eqmer2})
and which is also strongly related to the footprint bound from
Gr\"{o}bner basis theory.\\

 Given a zero dimensional ideal of a
multivariate polynomial ring, and a fixed 
monomial ordering, the well-known footprint bound states that the size of the
corresponding variety is at most equal to the number of monomials that
can not be found as leading monomial of any polynomial in the
ideal (if moreover the ideal is radical, then equality holds). 
More details on the footprint bound can be found in
\cite{clo,onorin,gh} -- in particular see \cite[Pro.\ 4, Sec.\ 5.3]{clo}.
We
 have the following easy corollary. \\ 

 \begin{Corollary}\label{footprspecial}
Given a polynomial $F(X_1, \ldots , X_m) \in {\mathbb{F}}[X_1,
\ldots , X_m]$, and a monomial ordering, let $X_1^{i_1}\cdots
X_m^{i_m}$ be the leading monomial of $F$, and assume $i_1 < s_1,
\ldots , i_m < s_m$. The number of elements in
$S_1 \times \cdots \times S_m$ that are zeros of $F$ is at most
equal to 
\begin{equation} 
s_1 \cdots s_m-(s_1-i_1)(s_2-i_2)\cdots (s_m-i_m).\label{eqthisisfoot}
\end{equation}
\end{Corollary}

\startprf
The set of zeros of $F$ from $S_1 \times \cdots \times S_m$ equals the
variety of the ideal $\langle F, G_1, \ldots , G_m\rangle$ where
$G_i=\prod_{u=1}^{s_i}(X_i-\alpha_u^{(i)})$. Here, we used the
notation introduced prior to Proposition~\ref{pronuogsaa}. The above ideal clearly is
zero-dimensional. In fact, the monomials that are not leading monomial
of any polynomial in the ideal must belong to the set
$$\{X_1^{j_1} \cdots X_m^{j_m} \mid j_1 < s_1, \ldots, j_m < s_m,
X_1^{j_1} \cdots X_m^{j_m}{\mbox{ is not divisible by }} X_1^{i_1}
\cdots X_m^{i_m} \},$$
the size of which equals (\ref{eqthisisfoot}). 
The result now follows from the footprint bound. \qed\\

The above corollary and (\ref{eqmer2}) are clearly related as
(\ref{eqthisisfoot}) equals the right side of (\ref{eqmer2}) for
$m=2$. Similarly, (\ref{eqbett}) equals the assumption in the corollary. Observe, however, that in (\ref{eqmer2}), and in this paper in
general, we always assume that the monomial ordering is the
lexicographic ordering described in Theorem~\ref{prorec}. The
master theorem of the present paper is the following result where
(\ref{eqendnuenstjerne}) is the generalization of (\ref{eqmer2}) to
more variables and
where the mentioned Condition A is the generalization
of~(\ref{eqbett}). Recall that $D(i_1, \ldots , i_m,r,s_1, \ldots ,
s_m)$ serves as an upper bound on the number of zeros of multiplicity
at least $r$ for polynomials with leading monomial being $X_1^{i_1}
\cdots X_m^{i_m}$ with respect to the lexicographic ordering. As a consequence the master theorem also can be viewed as a
generalization of Corollary~\ref{footprspecial}, when restricted to a lexicographic ordering. 

\begin{Theorem}\label{prosmall}
Assume that $(i_1, \ldots , i_m,r,s_1, \ldots , s_m)$ with $m\geq 2$ satisfies Condition
A in Definition~\ref{defconda} below. We have 
\begin{eqnarray}
D(i_1, \ldots ,i_m,r,s_1,\ldots , s_m) \leq s_1\cdots s_m-(s_1-\frac{i_1}{r})\cdots
(s_m-\frac{i_m}{r})\label{eqendnuenstjerne}
\end{eqnarray}
which is at most equal to $\min \{ (i_1s_2\cdots s_m+\cdots +s_1\cdots
s_{m-1}i_m)/r, s_1\cdots s_m\}$.
\end{Theorem}

We postpone the proof of Theorem~\ref{prosmall} till the end of the
section. \\

\begin{Definition}\label{defconda}
Let $m \geq 2$. We say that $(i_1, \ldots , i_m,r,s_1, \ldots , s_m)$
satisfies Condition A if the following hold
$$
\begin{array}{rl}
(A.1)& 0 \leq i_1\leq s_1, \ldots , 0 \leq i_{m-1} \leq s_{m-1}, 0
  \leq i_m < r s_m\\
(A.2)&s(s_1-\frac{i_1}{\ell}) \cdots (s_{m-2}-\frac{i_{m-2}}{\ell}) \leq \ell
  (s_1-\frac{i_1}{s})\cdots (s_{m-2}-\frac{i_{m-2}}{s})\\
& {\mbox{ \ for all \ }} 
  \ell=2, \ldots , r, s=1, \ldots \ell-1.\\
(A.3)&s(s_1-\frac{i_1}{r}) \cdots (s_{m-1}-\frac{i_{m-1}}{r}) \leq r
  (s_1-\frac{i_1}{s})\cdots (s_{m-1}-\frac{i_{m-1}}{s})\\
&{\mbox{ \  for all \ }} s=1, \ldots
  , r-1.
\end{array}
$$
\end{Definition}

We note that one could actually replace $\ell=2, \ldots , r$ in (A.2)
with the weaker $\ell =2, \ldots, r-1$ as the case $\ell=r$ follows
from (A.3).\\

Admittedly, the definition of the exponent being small (Condition A) is
rather technical. However: 
\begin{itemize}
\item If $(i_1, \ldots , i_m)$ is small
then all $(i_1^\prime, \ldots , i_m^\prime)$ with
$i_1^\prime \leq i_1, \ldots , i_m^\prime \leq i_m$ are also small
(Proposition~\ref{proalso}). Hence, it is enough to check if $(i_1,
\ldots , i_m)$ satisfies Condition A.
\item Condition A is satisfied when $i_t \leq s_t \min \left\{ \frac{\sqrt[m-1]{r}-1}{\sqrt[m-1]{r}-\frac{1}{r}},
\frac{\sqrt[m-2]{2}-1}{\sqrt[m-2]{2}-\frac{1}{2}}\right\}$, $t=1, \ldots
  ,m-1$, $i_m <r s_m$ (Theorem~\ref{themain}).
\item As already mentioned, Condition A and the master theorem reduces
  to well-known results when $r=1$ or when $m=2$ (see Remark~\ref{remsmart}
  for the details). 
\item For arbitrary $m$ but 
$r=2$ and $s_1 = \cdots =s_m$, 
Condition A reduces
to a simple expression (Proposition~\ref{proI} and Example~\ref{exsmalll}).
\end{itemize}

\begin{Proposition}\label{proalso}
If $(i_1, \ldots , i_m,r,s_1, \ldots , s_m)$ satisfies Condition
A then for all $i_1^\prime, \ldots ,i_m^\prime$ with $0 \leq
i_1^\prime \leq i_1, \ldots , 0 \leq i_m^\prime \leq i_m$ also
$(i_1^\prime , \ldots , i_m^\prime , r, s_1, \ldots , s_m)$ satisfies
Condition A.    
\end{Proposition}
\startprf
It is enough to show that
\begin{equation}
\frac{s_t-\frac{i_t}{s}}{s_t-\frac{i_t}{\ell}} \leq \frac{s_t-\frac{ai_t}{s}}{s_t-\frac{ai_t}{\ell}} \label{eqa}
\end{equation}
holds for all rational numbers $a$ and integers $t$  with $0 < a < 1$
and $1 \leq t \leq m-1$. But~(\ref{eqa}) is equivalent to
$(1-a)(\ell-s)\geq 0$
which is a valid inequality when $\ell >s$. \qed
\ \\

We now give the most important theorem of the paper.

\begin{Theorem}\label{themain}
If $i_m < r s_m$ and if for $t=1, \ldots , m-1$
$$i_t \leq s_t \min \left\{ \frac{\sqrt[m-1]{r}-1}{\sqrt[m-1]{r}-\frac{1}{r}},
\frac{\sqrt[m-2]{2}-1}{\sqrt[m-2]{2}-\frac{1}{2}}\right\}$$
then $D(i_1, \ldots  i_m,r,s_1, \ldots  s_m) \leq s_1 \cdots
s_m-(s_1-\frac{i_1}{r}) \cdots (s_m-\frac{i_m}{r})$.
\end{Theorem}
\startprf
The idea behind Theorem~\ref{themain} is to choose $i_t$, $t=1, \ldots
, m-1$ such that
\begin{equation}
\sqrt[m-1]{s}(s_t-\frac{i_t}{r}) \leq \sqrt[m-1]{r}(s_t-\frac{i_t}{s}),
{\mbox{ \ \ for }} s=1, \ldots , r-1, \label{TREKANT1}
\end{equation}
and such that 
\begin{equation}
\sqrt[m-2]{s}(s_t-\frac{i_t}{\ell}) \leq \sqrt[m-2]{\ell}(s_t-\frac{i_t}{s}),
{\mbox{ \ \ for }} \ell=2,\ldots , r, {\mbox{ / }}s=1, \ldots , \ell-1. \label{TREKANT2}
\end{equation}
The first set of inequalities guarantees (A.3) and the second set
guarantees (A.2). Now (\ref{TREKANT1}) and (\ref{TREKANT2}),
respectively, translates to
\begin{equation}
\frac{i_t}{s_t}  \leq 
\frac{\sqrt[m-1]{r}-\sqrt[m-1]{s}}{\frac{\sqrt[m-1]{r}}{s}-\frac{\sqrt[m-1]{s}}{r}}, \label{nabel1}
\end{equation}
\begin{equation}
\frac{i_t}{s_t}  \leq 
\frac{\sqrt[m-2]{\ell}-\sqrt[m-2]{s}}{\frac{\sqrt[m-2]{\ell}}{s}-\frac{\sqrt[m-2]{s}}{\ell}}, \label{nabel2}
\end{equation}
respectively, and then also (A.1) is clearly satisfied. We shall show
that the right side of (\ref{nabel1}) is smallest possible when $s=1$,
in which case it equals
$(\sqrt[m-1]{r}-1)/(\sqrt[m-1]{r}-1/r)$. And we shall show that the
right side of (\ref{nabel2}) is smallest possible when $\ell=2, s=1$,
in which case it equals $(\sqrt[m-2]{2}-1)/(\sqrt[m-2]{2}-1/2)$.\\
We first consider (\ref{nabel1}) where we substitute $S=\sqrt[m-1]{s}$
and $R=\sqrt[m-1]{r}$ to obtain
$$\frac{i_t}{s_t} \leq \frac{R^mS^{m-1}-R^{m-1}S^m}{R^m-S^m}.$$
We want to demonstrate that the right side is minimal on $[ 1, R [$
when $S=1$. The derivative is
$$\frac{(m-1)R^{2m}S^{m-2}+R^mS^{2m-2}-mR^{2m-1}S^{m-1}}{(R^m-S^m)^2}.$$
Hence, it suffices to show that the numerator is always positive on
$] 0,R[$. Writing $S=Ra$ with $a \in ] 0, 1[$ the condition that
the numerator should be positive becomes $m-1+a^m-ma >0$. Plugging in
$a=1$, equality holds. Therefore the result follows from the fact that
the derivative of $m-1+a^m-ma$ is negative on $]0,1[$.\\
The above proof not only shows that the minimum of the right side of
(\ref{nabel1}) is obtained for $s=1$. It also applies to demonstrate
that the minimum of the right side of (\ref{nabel2}) is attained in
one of the following cases $(\ell=2, s=1)$, $(\ell=3, s=1), \ldots ,
(\ell=r,s=1)$. We next substitute $m-2$ with $m$ on the right side of
(\ref{nabel2}) to obtain $(\ell^{1/m}-1)/(\ell^{1/m}-1/\ell)$. We want
to show that the minimal value for $\ell \in [2, \infty [$ is
attained when $\ell=2$. The derivative is
$$\frac{-\left(\ell^{1/m}m-\ell^{(m+1)/m}+\ell^{1/m}-m\right)}{m
  \left(\ell^{(m+1)/m}-1\right)^2}$$
where the denominator is always positive and the numerator is
positive for $\ell=0$. The result follows from the fact that 
$$\frac{d}{d\ell}\left(\ell^{1/m}m-\ell^{(m+1)/m}+\ell^{1/m}-m)\right)=\frac{(m+1)(\ell^{-(m-1)/m}-\ell^{1/m})}{m}$$
is negative on $]0,\infty[$. 
\qed

\begin{Remark}\label{remsmart}
If $r=1$ then (A.2) and (A.3) do not apply and therefore Condition A
reduces to $i_1 \leq s_1, \ldots , i_m \leq s_m$. Hence, in this case Theorem~\ref{prosmall} in combination with Theorem~\ref{prorec} reduce to
Theorem~\ref{footprspecial}.\\ 
For $m=2$ and $r$ arbitrary  condition (A.2) does not apply and
condition (A.3)
simplifies to
$$i_1 \leq \frac{rs}{r+s}s_1$$
for all integers $s$ with $1 \leq s <  r$. The minimal upper bound on
$i_1$ is attained for $s=1$. Hence, in case of two variables Condition A reads $i_1
\leq \frac{r}{r+1}s_1$, $i_2 < r s_2$. For $m=2$ and $r$ being arbitrary Theorem~\ref{prosmall} therefore equals (\ref{eqmer2}) and (\ref{eqbett}).
\end{Remark}

\begin{Proposition}\label{proI}
Assume $r=2$ and $s_1=\cdots =s_m =q$. Then Condition A simplifies to
$$\sum_{t=1}^{m-1} (-1)^{t+1} \frac{2^{t+1}-1}{2^t}\sum_{1 \leq j_1 <
  \cdots < j_t \leq m-1} (I_{j_1} \cdots I_{j_t}) \leq 1 {\mbox{ \ and \ }} I_m < 2$$
where $I_1=i_1/q, \ldots , I_m=i_m/q$.  
\end{Proposition}
\startprf
For $r=2$, the conditions (A.2), (A.3) become  
$$\big(s_1-\frac{i_1}{2}\big) \cdots\big(s_{m-1}-\frac{i_{m-1}}{2}\big) \leq
2\big(s_1-i_1\big) \cdots\big(s_{m-1}-i_{m-1}\big)$$
which is equivalent to
\begin{eqnarray}
&&\left( 1-\frac{I_1}{2}\right)\cdots \left(
  1-\frac{I_{m-1}}{2}\right)\leq 2(1-I_1)\cdots (1-I_{m-1}) \nonumber
  \\
&\Updownarrow \nonumber \\
&&1+\sum_{t=1}^{m-1}(-1)^t(\frac{1}{2})^t\sum_{1 \leq j_1 < \cdots < j_t
    \leq m-1} (I_{j_1} \cdots I_{j_t})\leq \nonumber \\
&&{\mbox{ \ \ \ \ \ \ }} 2 +
  2\sum_{t=1}^{m-1}(-1)^t\sum_{1 \leq j_1 < \cdots <j_t \leq m-1}(I_{j_1}
  \cdots I_{j_t}) \nonumber \\
&\Updownarrow \nonumber \\
&&\sum_{t=1}^{m-1} (-1)^{t+1} \frac{2^{t+1}-1}{2^t}\sum_{1 \leq j_1 <
  \cdots < j_t \leq m-1} (I_{j_1} \cdots I_{j_t}) \leq 1 \nonumber
\end{eqnarray}
and we are through. \qed

\begin{Example}\label{exsmalll}
Let the notation be as in Proposition~\ref{proI}. 
For $r=2$, $m=3$ and $s_1=s_2=s_3=q$ Condition A reads
$$\frac{3}{2}(I_1+I_2)-\frac{7}{4}I_1I_2 \leq 1, {\mbox{ \ \ }} I_3
< 2.$$ For $r=2$, $m=4$ and $s_1=s_2=s_3=s_4=q$ Condition A reads
  $$\frac{3}{2}(I_1+I_2+I_3)-\frac{7}{4}(I_1I_2+I_1I_3+I_2I_3)+\frac{15}{8}I_1I_2I_3 \leq 1, {\mbox{ \ \ }} I_4 < 2.$$
This is illustrated in Figure~\ref{figimp}. 
\begin{figure}
\begin{center}
\input{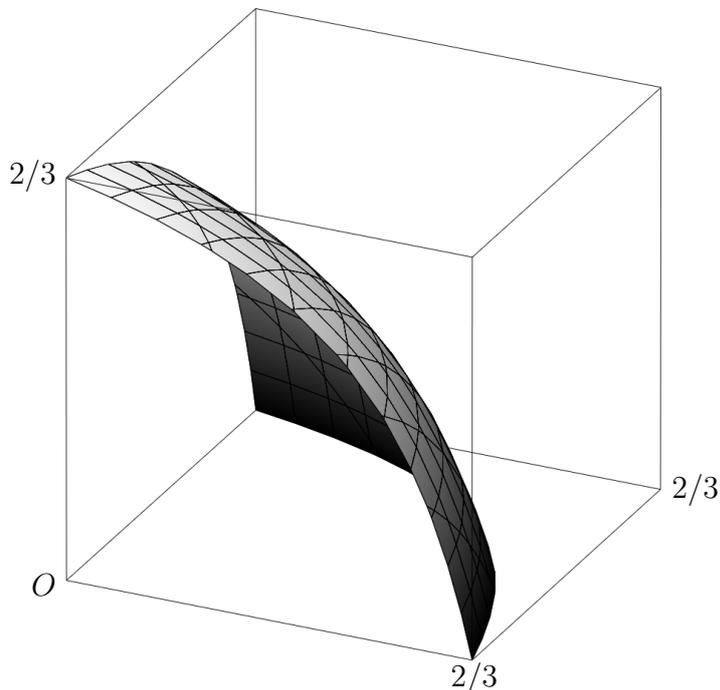}
\end{center}
\caption{The surface $\frac{3}{2}(I_1+I_2+I_3)-\frac{7}{4}(I_1I_2+I_1I_3+I_2I_3)+\frac{15}{8}I_1I_2I_3 = 1$}
\label{figimp}
\end{figure}
\end{Example}

From Proposition~\ref{proI} it is clear that in the case of $r=2$, for
Condition A to hold we must have $i_t \leq \frac{2}{3} s_t$, $t=1,
\ldots , m-1$. The general picture for $r$ arbitrary is described in
the following proposition.

\begin{Proposition}\label{lem}
Assume that $(i_1, \ldots , i_m,r,s_1, \ldots , s_m)$ with $m\geq 2$ satisfies Condition
A. If $r \geq 2$ then 
\begin{equation}
i_1 \leq \frac{r}{r+1}s_1, \ldots , i_{m-1} \leq \frac{r}{r+1}s_{m-1}.\label{eqeqeq}\end{equation}
\end{Proposition}
\startprf
Follows from (A.3), the last part of Remark~\ref{remsmart}, and the fact that 
$$s_t-\frac{i_t}{\ell}\geq s_t-\frac{i_t}{s}$$
holds for $t=1, \ldots , m-1$. \qed

\ \\

\startprfoftheorem
Let $(i_1, \ldots ,i_m,r,s_1, \ldots ,s_m)$ with $m\geq 2$ be
such that Condition A holds. We give an induction proof
that 
\begin{equation}
\begin{array}{r}
D(i_1, \ldots , i_t,l,s_1,\ldots ,s_t)  \leq s_1\cdots s_t-(s_1-\frac{i_1}{\ell})
\cdots (s_t-\frac{i_{t}}{\ell})\\
 {\mbox{ \ for all \ }}1 \leq t < m , 1 \leq
  \ell\leq r.
\end{array}
\label{eqminus1}
\end{equation} 
For $t=1$ the
result is clear. Let $1 < t<m$ and assume the result holds when $t$ is
substituted with $t-1$. According to Definition~\ref{defD} we have
\begin{multline*}
D(i_1, \ldots , i_t,l,s_1, \ldots ,s_t)= 
\\
\begin{split}
\max_{(u_1, \ldots  ,u_l)\in A(i_t,\ell,s_t) }\bigg\{& (s_t-u_1-\cdots -u_\ell)D(i_1,\ldots ,i_{t-1},\ell,s_1,
\ldots ,s_{t-1})\\
&+u_1D(i_1, \ldots , i_{t-1},\ell-1,s_1, \ldots ,s_{t-1})+\cdots 
\\
&+u_{\ell-1}D(i_1, \ldots ,i_{t-1},1,s_1, \ldots ,
s_{t-1})+u_\ell s_1\cdots s_{t-1} \bigg\} 
\end{split}
\end{multline*}
where
\begin{eqnarray}
A(i_t,\ell,s_t)&=&\{(u_1, \ldots ,u_\ell) \in {\mathbb{N}}_0^\ell \mid
u_1+ \cdots +u_\ell \leq s_t, {\mbox{ \ }}u_1+2u_2+\cdots + \ell u_\ell \leq i_t\}\nonumber
\end{eqnarray}
follows from~ Definition~\ref{defD}. By the above assumptions this implies that
\begin{multline}
D(i_1, \ldots ,i_t,\ell,s_1, \ldots ,s_t) \leq
\\
\shoveleft
\max_{(u_1, \ldots  ,u_\ell)\in B(i_t,\ell,s_t) } \bigg\{ s_t\big(
s_1 \cdots s_{t-1}-(s_1-\frac{i_1}{\ell})\cdots(s_{t-1}-\frac{i_{t-1}}{\ell})\big)\\
\begin{split}
&+u_1\big((s_1-\frac{i_1}{\ell})\cdots(s_{t-1}-\frac{i_{t-1}}{\ell})-(s_1-\frac{i_1}{\ell-1})\cdots(s_{t-1}-\frac{i_{t-1}}{\ell-1}) \big)\\
&+ \cdots
\\
&+u_{\ell-1}\big(
(s_1-\frac{i_1}{\ell})\cdots(s_{t-1}-\frac{i_{t-1}}{\ell})-(s_1-\frac{i_1}{1})\cdots(s_{t-1}-\frac{i_{t-1}}{1})
\big)\\
&+u_{\ell}\big((s_1-\frac{i_1}{\ell})\cdots(s_{t-1}-\frac{i_{t-1}}{\ell})\big) \bigg\}\label{eqznabel}
\end{split}
\end{multline}
where 
\begin{eqnarray}B(i_{t},\ell,s_t)&=&\{(u_1, \ldots , u_\ell) \in {\mathbb{Q}}^\ell \mid 0
\leq u_1, \ldots , u_\ell, {\mbox{ \ }} u_1 + \cdots +u_\ell \leq s_t,
                                     \nonumber \\
&&{\mbox{ \ \hspace{4cm} and \ }} u_1+2u_2+\cdots +\ell u_\ell \leq i_t\}.\nonumber
\end{eqnarray}
We have $t<m $ and therefore condition  (A.2) applies. We note that 
$$s(s_1-\frac{i_1}{\ell}) \cdots (s_{t-1}-\frac{i_{t-1}}{\ell})
\leq \ell (s_1-\frac{i_1}{s})\cdots (s_{t-1}-\frac{i_{t-1}}{s})$$
for $s=1, \ldots , \ell-1$ is equivalent to 
$$(\ell-s)(s_1-\frac{i_1}{\ell}) \cdots (s_{t-1}-\frac{i_{t-1}}{\ell})
\leq \ell (s_1-\frac{i_1}{\ell-s})\cdots (s_{t-1}-\frac{i_{t-1}}{\ell-s})$$
for $s=1, \ldots , \ell-1$ which again is equivalent to
$$\ell\big( (s_1-\frac{i_1}{\ell}) \cdots
(s_{t-1}-\frac{i_{t-1}}{\ell})-(s_1-\frac{i_1}{\ell-s}) \cdots
(s_{t-1}-\frac{i_{t-1}}{\ell-s}) \big) 
\leq s (s_1-\frac{i_1}{\ell})\cdots (s_{t-1}-\frac{i_{t-1}}{\ell})$$
for $s=1, \ldots , \ell-1$. Therefore the maximal value
of~(\ref{eqznabel}) is attained for $u_1=\cdots =u_{\ell-1}=0$ and
$u_\ell=\frac{i_t}{\ell}$. This concludes the induction proof of (\ref{eqminus1}).\\
To show (\ref{eqendnuenstjerne}) we apply similar arguments
to the case
$t=m$
but use condition (A.3) rather than condition (A.2). \\
Finally we address the last part of Theorem~\ref{prosmall}. It is
clear that the right side of~(\ref{eqendnuenstjerne}) is smaller than or equal to $s_1\cdots s_m$. To
see that it is also smaller than or equal to
\begin{equation}
\sum_{t=1}^m \big( (\prod_{\begin{array}{c}j=1,\ldots , m\\j \neq t
\end{array}}s_j)\frac{i_t}{r}\big) \label{eqsumprod}
\end{equation}
we start by observing that 
$$\big( \prod_{\begin{array}{c}j=1, \ldots , m\\j \neq
    t\end{array}}s_j \big) \frac{i_t}{r}$$
 equals the volume of 
\begin{eqnarray}
N(t,\frac{i_t}{r})&=&\{(a_1, \ldots , a_m) \in {\mathbb{R}}_0^m \mid 0
\leq a_t < \frac{i_t}{r}, 0 \leq a_j \leq s_j \nonumber \\
&& {\mbox{ \ \ \ \ \ \ \ \ \ \ \ \ \ \ \ \ \ \ \ \ \ \ \ \ \ \ \ \  \ \ \  for }} j \in \{1,
\ldots , m\}\backslash \{ t \} \}. \nonumber
\end{eqnarray}
The sum of volumes of $N(t,\frac{i_t}{r})$, $t=1, \ldots m$ is larger
than or equal to the volume of
\begin{eqnarray}
\cup_{t=1}^m N(t, \frac{i_t}{r})&=&\{(a_1, \ldots ,a_m) \in
    {\mathbb{R}}_0^m \mid 0 \leq a_t \leq s_t {\mbox{ for }} t=1, \ldots
    , m \nonumber \\
&&{\mbox{ \ \ \ \ \  \ \ \ \ \ \ \ \ \ \ \ \ \ \ \ \ \ \ \ \ \ \ \ \  and not all $j$ satisfy }} \frac{i_j}{r} \leq a_j\} \nonumber
\end{eqnarray}
which equals the right side of~(\ref{eqendnuenstjerne}). \qed

\section{Concluding remarks}
 The results in this paper use the lexicographic
ordering. We pose it as a research problem to investigate if some of
them hold for arbitrary monomial orderings. 

\section*{Acknowledgments}
This work was supported by the Danish Council for Independent
Research (grant no.\ DFF-4002-00367) and by the Danish National
Research Foundation and the National Natural Science
Foundation of China (Grant No.\ 11061130539 -- the Danish-Chinese
Center for Applications of Algebraic Geometry in Coding Theory and
Cryptography).

\end{document}